\documentclass{article}
\usepackage{amsmath}
\usepackage{latexsym}
\usepackage{amssymb}
\usepackage[dvips]{graphicx}
\newtheorem{lem}{Lemma}
\newtheorem{thm}{Theorem}
\newtheorem{cor}{Corollary}
\newtheorem{prop}{Proposition}

\def\QED{\mbox{\rule[0pt]{1.5ex}{1.5ex}}}

\def\endproof{\hspace*{\fill}~\QED\par\endtrivlist\unskip}

\def\inner{\mathop{\rm int}\nolimits}

\def\real{\mathbb{R}}

\def\Label#1{\label{#1}\ [\ #1\ ]\ }
\def\Label{\label}
\def\Bibitem#1{\bibitem{#1}\ [\ #1\ ]\ }
\def\Bibitem{\bibitem}
\begin{document}
\begin{center}
\bf\Large Two non-regular extensions of the large deviation bound
\end{center}
\begin{center}
\bf\large
Masahito Hayashi\footnote{e-mail masahito@brain.riken.go.jp}\\
Laboratory for Mathematical Neuroscience, 
Brain Science Institute, RIKEN
\footnote{2-1 Hirosawa, Wako, Saitama, 351-0198, Japan}
\end{center}
\begin{abstract}
We formulate two types of extension of the large deviation theory initiated by Bahadur in a non-regular setting. One can be regarded as a bound of the point estimation, the other can be regarded as the limit of a bound of the interval estimation. Both coincide in the regular case, but do not necessarily coincide in a non-regular case. Using the limits of relative R\'{e}nyi entropies, we derive their upper bounds and give a necessary and sufficient condition for the coincidence of the two upper bounds. We also show the attainability of these two bounds in several non-regular location shift families.
\end{abstract}
\section{Introduction}\Label{s1}
As discussed by Bahadur \cite{Ba60,Ba67,Ba71}, Fisher information characterizes the limit of the decreasing rate of the tail probability of the optimal estimator. However, 
when the Kullback-Leibler (KL) divergence diverges or Fisher information cannot be defined, this cannot be applied. As an 
alternative information quantity between two probabilities, we can use the relative R\'{e}nyi entropies $I^s(p\|q)
:= -\log \int p^s(\omega) q^{1-s}(\omega) \,d \omega ~
(0 \,< s \,< 1)$, which play an important role in simple hypothesis testing. (Hoeffding \cite{Hoef}, 
Chernoff \cite{Cher})
In general, for a probability distribution 
family $\{p_{\theta}|\theta \subset \real\}$, the relative R\'{e}nyi entropies $I^s(p_{\theta}\|p_{\theta+\epsilon})$ tend to $0$ when $ \epsilon$ goes to $0$, but their order is not necessarily $\epsilon^2$. In this paper, we treat a large deviation theory, which can be applied to such a case. The importance of the relative R\'{e}nyi entropies $I^s(p_{\theta}\|p_{\theta+\epsilon})$ was pointed out by Akahira \cite{Ak} 
from the viewpoint of information loss in a non-regular family. Although the limit distribution of the maximum likelihood estimator (MLE) in a non-regular family has been extensively discussed \cite{Smith,Cheng,AT2}, large deviation theory has not been discussed sufficiently in this respect.

In a large deviation theory of the parameter estimation, we usually focus on the first exponential rate
\begin{align}
\beta(\vec{T},\theta,\epsilon):=
\liminf \frac{-1}{n} \log p^n_\theta
\{ | T_n - \theta | \,> \epsilon\} \Label{12}
\end{align}
for a sequence of estimators $\vec{T}=\{T_n\}$, which is simply called an estimator in the following. Of course, when $\beta(\vec{T},\theta,\epsilon)$ is large, the estimator $\vec{T}$ is better. We discuss its maximization at the limit $\epsilon \to 0$.
As is explained in Section \ref{s2}, Bahadur focused on the bound $\alpha(\theta)$ of the first exponential rate at the limit $\epsilon \to 0$, whose definition is precisely given in Section \ref{s2}.
He derived its upper bound from the viewpoint of Stein's lemma in simple hypothesis testing.

The main purpose of this paper is to extend the bound $\alpha(\theta)$ and derive its upper bound in a more general setting. Indeed, depending on how the limit $\epsilon \to 0$ is obtained, we can extend the bound $\alpha(\theta)$ in two ways. This difference can be regarded as the difference between the point estimation and the limit of the interval estimation. Bahadur's original theory concerned the point estimation, but it can also be applied to the limit of the interval estimation, as is explained in Section \ref{s2}.
Therefore, from the two ways of extending the bound $\alpha(\theta)$, we can define two generalizations $\alpha_1(\theta)$ and $\alpha_2(\theta)$ of Bahadur's bound $\alpha(\theta)$. In Section \ref{s3}, we give the respective upper bounds $\overline{\alpha}_1(\theta)$ and $\overline{\alpha}_2(\theta)$, and check that our result contains Bahadur's results as a special case in which the two upper bounds $\overline{\alpha}_1(\theta)$ and $\overline{\alpha}_2(\theta)$ coincide. Proofs of the two inequalities $\alpha_1(\theta) \le \overline{\alpha}_1(\theta)$ and $\alpha_2(\theta) \le \overline{\alpha}_2(\theta)$ are given in Section \ref{s4}. 

In Sections \ref{s5} and \ref{s6}, 
we also discuss the attainability of upper bounds $\overline{\alpha}_1(\theta)$ and $\overline{\alpha}_2(\theta)$, and calculate them for location shift families. In Section \ref{s5}, 
we derive several formulae for the first exponential rates of the maximum likelihood estimator (MLE) and some estimators consisting of order statistics under certain assumptions. 
In Section \ref{s6}, using a formula for relative R\'{e}nyi entropies given by Hayashi \cite{Haya2}, we calculate the two upper bounds 
$\overline{\alpha}_1(\theta)$ and $\overline{\alpha}_2(\theta)$, and derive a sufficient and necessary condition for their coincidence. Using formulae given in Section \ref{s5}, we check that these bounds are attainable in some special cases. Our examples in Section \ref{s6} contain location shift families generated by a Weibull distribution, gamma distribution, beta distribution, and uniform distribution.

\section{Bahadur theory}\Label{s2}
In this section, we begin by summarizing the results reported by Bahadur \cite{Ba60,Ba67,Ba71}, who discussed the decreasing rate of the tail probability in the estimation for a distribution family. Given $n$-i.i.d. data 
$\omega_1, \ldots, \omega_n$, the first exponential rate $\beta(\vec{T},\theta,\epsilon)$ of the estimator 
$\vec{T}=\{T_n\}$ is written as 
\begin{align*} \beta(\vec{T},\theta,\epsilon)= \min \{ \beta^+(\vec{T},\theta,\epsilon), \beta^-(\vec{T},\theta,\epsilon)\},
\end{align*}
where the exponential rates of half-side error
probabilities are given by
\begin{align*}
\beta^+(\vec{T},\theta,\epsilon) &:=
\liminf \frac{-1}{n} \log p^n_\theta
\{ T_n \,> \theta + \epsilon\} \\
\beta^-(\vec{T},\theta,\epsilon) &:=
\liminf \frac{-1}{n} \log p^n_\theta
\{ T_n \,< \theta - \epsilon\}.
\end{align*}
When an estimator $\vec{T}=\{T_n\}$ satisfies the weak consistency
\begin{align*}
 p^n_\theta
\{ | T_n - \theta | \,> \epsilon\} \to 0 
\quad \forall \epsilon \,> 0, \quad
\forall \theta \in \Theta,
\end{align*}
using the monotonicity of KL-divergence, we can prove the inequality
\begin{align}
\beta(\vec{T},\theta,\epsilon)
\le \min\{ D(p_{\theta+\epsilon} \|p_\theta),
D(p_{\theta-\epsilon} \|p_\theta)\}. \Label{Baha4}
\end{align}
Note that if, and only if, the family is exponential, there exists an estimator attaining the equality
(\ref{Baha4}) at 
$\forall\theta\in \Theta, \forall\epsilon \,> 0$.
Therefore, for a general family, it is difficult to optimize the first exponential rate $\beta(\vec{T},\theta,\epsilon)$.

We usually introduce the second exponential rate as another optimized value:
\begin{align}
\alpha(\vec{T},\theta)
:= 
\lim _{\epsilon \to + 0}
\frac{1}{\epsilon^2}
\beta(\vec{T},\theta,\epsilon). \Label{2.5}
\end{align}
In this case, when the Fisher information $J_\theta$ satisfies the condition
\begin{align*}
J_\theta:= \int l_\theta (\omega)^2 p_\theta(\,d \omega)
=\lim _{\epsilon \to 0}
\frac{2}{\epsilon^2}
D(p_{\theta+\epsilon} \|p_\theta),
\quad l_\theta (\omega):=
\frac{\,d }{\,d \theta}
\log p_{\theta} (\omega),
\end{align*}
the inequality
\begin{align}
\alpha(\vec{T},\theta)
\le
\frac{1}{2} J_\theta
\Label{Baha1}
\end{align}
holds.
Moreover, as was proven by Fu \cite{Fu73}, if the family satisfies the concavity of the logarithmic derivative $l_\theta(\omega)$ for $\theta$ and some other conditions, the MLE $\vec{\theta_{ML}}$ attains the equality of (\ref{Baha1}).
These facts are summarized in the two inequalities
\begin{align}
\alpha(\theta) :=
\sup_{\vec{T}:{\rm WC}}
\alpha(\vec{T},\theta) 
&=
\frac{1}{2} J_\theta \Label{Baha2} \\
\lim _{\epsilon \to + 0}
\frac{1}{\epsilon^2}
\sup_{\vec{T}:{\rm WC}}
\beta(\vec{T},\theta,\epsilon)
&=
\frac{1}{2} J_\theta \Label{Baha3}.
\end{align}
Inequality (\ref{Baha2}) can be regarded as the bound of the point estimation, while inequality (\ref{Baha3}) can be regarded as the limit of the bound of the interval estimation because $\sup_{T:WC}\beta(\vec{T},\theta,\epsilon)$ corresponds to the bound of the interval estimation with width $2 \epsilon$ of the confidence interval.
\section{Main results}\Label{s3}
In this paper, the relative R\'{e}nyi entropies $I^s(p_\theta\|p_{\theta+\epsilon})$ substitute for the KL divergence. Note that the order of $I^s(p_\theta\|p_{\theta+\epsilon})$ is not necessarily $\epsilon^2$ at the limit $\epsilon \to 0$.
However, its order is independent of the parameter $s$, as is guaranteed by the inequalities
\begin{align*}
2 \min \{ s, 1-s\} I^{\frac{1}{2}}
(p_\theta\|p_{\theta+\epsilon})
\le I^s(p_\theta\|p_{\theta+\epsilon})
\le
2 \max \{ s, 1-s\} I^{\frac{1}{2}}
(p_\theta\|p_{\theta+\epsilon}),
\end{align*}
which are proven in Lemma \ref{L38} of the Appendix. In several cases, the order of the first exponential rate
$\beta(\vec{T},\theta,\epsilon)$ coincides with the order of the relative R\'{e}nyi entropies $I^s(p_\theta\|p_{\theta+\epsilon})$.
In the following, we use a strictly monotonically decreasing function $g(x)$ such that $I^s(p_\theta\|p_{\theta+\epsilon})\cong O(g(\epsilon)) $ and $g(0)=0$. 

Following equations (\ref{Baha2}) and (\ref{Baha3}), 
we define two extensions of Bahadur's bound $\alpha(\theta)$ as
\begin{align}
\alpha_1(\theta)&:=
\limsup_{\epsilon \to + 0}
\frac{1}{g(\epsilon)}
\sup_{T}
\inf_{\theta-\epsilon \le \theta'\le \theta+\epsilon}
\beta(\vec{T},\theta',\epsilon) \Label{hi1} \\
\alpha_2(\theta)&:=
\sup_{T}
\alpha_2(\vec{T},\theta),
\Label{hi2}
\end{align}
where
\begin{align}
\alpha_2(\vec{T},\theta):= \liminf_{\epsilon \to + 0}
\frac{1}{g(\epsilon)}
\inf_{\theta-\epsilon \le \theta'\le \theta+\epsilon}
\beta(\vec{T},\theta',\epsilon) .\Label{8}
\end{align}
Note that we take infimum
$\inf_{\theta-\epsilon \le \theta'\le \theta+\epsilon}$
into account in (\ref{8}), unlike (\ref{2.5}).
As was pointed out by Ibragimov and Has'minskii
\cite{IH},
when KL-divergence is infinite, there exists a consistent super efficient estimator $\vec{T}$ such that $\beta(\vec{T},\theta,\epsilon)$ and $\lim _{\epsilon \to + 0}
\frac{1}{g(\epsilon)} \beta(\vec{T},\theta,\epsilon)$ is infinite at one point $\theta$.
Therefore, we need to take the infimum 
$\inf_{\theta-\epsilon \le \theta'\le \theta+\epsilon}$
into account. In this situation, we do not need to 
limit estimators to weakly consistent ones. As is proven in the next section, we can obtain the following theorems. 
\begin{thm}\Label{thm1}
When the convergence $\lim_{\epsilon \to 0 } \frac{I^s(p_{\theta-\epsilon/2}\| p_{\theta+\epsilon/2})} {g(\epsilon)}$ 
is uniform for $0\,< s \,< 1$, the inequality
\begin{align}
\alpha_1(\theta) \le 
\overline{\alpha}_1(\theta):= 2^{\kappa}
\sup_{0 \,< s \,< 1}I^s_{g,\theta} \Label{7}
\end{align}
holds, where $\kappa$ and 
$I^s_{g,\theta}$ are defined by
\begin{align}
I^s_{g,\theta} &:= \lim_{\epsilon \to +0 } 
\frac{I^s(p_{\theta-\epsilon/2}\| 
p_{\theta+\epsilon/2})}
{g(\epsilon)} \quad 1 \ge s \ge 0 \Label{uni} \\
x^\kappa &= \lim_{\epsilon \to +0} \frac{g(x\epsilon)}
{g(\epsilon)} . \nonumber
\end{align}
Lemma \ref{ap1} proven in Appendix \ref{A}, 
guarantees the existence of such a real number $\kappa$.
\end{thm}
Note that the function $s \to I^s_g$ is concave and continuous because the function $s \to I^s(p_{\theta-\frac{1}{2} \epsilon}\|p_{\theta+\frac{1}{2}\epsilon})$ is concave and continuous.
Therefore, when $I_{g,\theta}^s=I_{g,\theta}^{1-s}$, we have
\begin{align}
\overline{\alpha}_1(\theta)= 2^\kappa
I^{\frac{1}{2}}_{g,\theta}.
\end{align}
\begin{thm}\Label{thm2}
If the convergence $\lim_{\epsilon \to 0 } \frac{I^s(p_{\theta-\epsilon/2}\| p_{\theta+\epsilon/2})} {g(\epsilon)}$ is uniform for $ s \in (0,1)$ and $\theta \in K$ for any compact set $K \subset \real$, the inequality
\begin{align}
\alpha_2(\theta)
\le
\overline{\alpha}_2(\theta) 
:=
\left\{
\begin{array}{cc}
\sup_{0\,< s \,< 1}
\frac{I^s_{g,\theta}}{s(1-s)}
\left( s^{\frac{1}{\kappa-1}} + (1-s)^{\frac{1}{\kappa-1}} 
\right)^{\kappa-1}& \hbox{ if }
\kappa \,< 1 \\
2 I^{\frac{1}{2}}_{g,\theta} & \hbox{ if }
\kappa = 1 \\
\inf_{0 \,< s \,< 1}
\frac{I^s_{g,\theta}}{s(1-s)}
\left( s^{\frac{1}{\kappa-1}} + (1-s)^{\frac{1}{\kappa-1}} 
\right)^{\kappa-1}& \hbox{ if }
\kappa \,> 1 
\end{array}
\right. 
\Label{9}
\end{align}
holds.
\end{thm}
In our proofs of these theorems, Chernoff's formula and Hoeffding's formula in simple hypothesis testing play important roles.

As was proven by Akahira \cite{Ak}, under some regularity conditions for a distribution family, the equation
\begin{align}
\lim_{\epsilon \to +0}
\frac{I^s(p_\theta\|p_{\theta+\epsilon})}
{\epsilon^2}
= \frac{J_\theta s(1-s)}{2}, \Label{38}
\end{align}
holds.
When we choose $g(x)=x^2$, we have $\kappa=2, I^s_{g,\theta} = \frac{1}{2}J_{\theta}s(1-s)$, and the relations
\begin{align*}
\overline{\alpha}_1(\theta) &= 
4 \max_{0\le s \le 1}
I^s_{g,\theta}
=\frac{1}{2} J_\theta \\
\overline{\alpha}_2(\theta) &= 
\min_{0\le s \le 1}
\frac{I^s_{g,\theta}}{s(1-s)}
=\frac{1}{2} J_\theta 
\end{align*}
hold.
In particular, if the distribution family satisfies the concavity of the logarithmic derivative 
$l_\theta(\omega)$ for $\theta$ and some other conditions, the bound $\frac{1}{2} J_\theta $ is attained by the MLE. Thus, 
the relations $\alpha_1(\theta) = \alpha_2(\theta) = \overline{\alpha}_1(\theta) = \overline{\alpha}_2(\theta) 
= \frac{1}{2} J_\theta$ hold.

As a relation between two bounds $\overline{\alpha}_1
(\theta) $ and $\overline{\alpha}_2(\theta)$, we can 
prove the following theorem.
\begin{thm}\Label{thm3}
The inequality
\begin{align}
\overline{\alpha}_1(\theta) 
\ge 
\overline{\alpha}_2(\theta)  \Label{14}
\end{align}
holds, and (\ref{14}) holds as an equality 
if and only if the equations 
\begin{align}
\overline{\alpha}_1(\theta) 
&=
2^\kappa I_{g,\theta}^{\frac{1}{2}} \Label{16.3} \\
2^\kappa I_{g,\theta}^{\frac{1}{2}}
&=
\overline{\alpha}_2(\theta) \Label{15}
\end{align}
hold.
When $\kappa \le 1$, 
(\ref{14}) holds as an equality 
if and only if equation (\ref{16.3}) holds. 
\end{thm}
When $I^s_{g,\theta}$ is differentiable, condition (\ref{16.3}) is equivalent to $ \left. \frac{\,d }{\,d s} {I}_{g,\theta}^{s}
\right|_{s=\frac{1}{2}} =0$.

\section{Proofs of main results}\Label{s4}
In our proofs of Theorems \ref{thm1} and \ref{thm2}, 
Chernoff's formula and Hoeffding's formula in simple 
hypothesis testing are essential, and are summarized 
as follows. 
Let the probability $p$ on $\Omega$ be the null 
hypothesis and $q$ be the alternative hypothesis. 
When we discuss 
a hypothesis testing problem concerning $n$-i.i.d. 
data, we call a sequence $\vec{A}=\{A_n\}$ a test, 
where $A_n$ is an acceptance region, which is a subset 
of $\Omega^n$. The first error probability $e_1(A_n)$ 
and the second error probability $e_2(A_n)$ are defined as
\begin{align*}
e_1(A_n):= 1- p^n(A_n) , \quad
e_2(A_n):= q^n(A_n),
\end{align*}
and their exponents are given by
\begin{align*}
e^*_1(\vec{A}):= \liminf_{n\to \infty}
\frac{-1}{n} \log e_1(A_n) \\
e^*_2(\vec{A}):= \liminf_{n\to \infty}
\frac{-1}{n} \log e_2(A_n) .
\end{align*}
Chernoff \cite{Cher} evaluated the exponent of the sum of the two errors as 
\begin{align}
\sup_{\vec{A}}
\lim_{n \to \infty} \frac{-1}{n}
\log
(e_1(A_n)+e_2(A_n))
=
\sup_{\vec{A}}
\min\{ e^*_1(\vec{A}),
e^*_2(\vec{A})\}
=\sup_{0\,< s \,< 1} I^s(p\|q),  \Label{11}
\end{align}
which is essential for our proof of Theorem \ref{thm1}.
This bound is achieved by both of the likelihood tests $\{ \vec{\omega}_n \in \Omega^n| p^n(\vec{\omega}_n) \ge q^n(\vec{\omega}_n)\}$ and $\{ \vec{\omega}_n \in \Omega^n| p^n(\vec{\omega}_n) \,> q^n(\vec{\omega}_n)\}$.

Hoeffding proved another formula for simple hypothesis testing \cite{Hoef}:
\begin{align}
\sup_{\vec{A}:e^*_1(\vec{A})\ge r }
e^*_2(\vec{A}) = \sup_{0 \,< s \,< 1}
\frac{-sr + I^s(p\|q)}{1-s}.
 \Label{Hoe}
\end{align}
This formula is essential for our proof of Theorem \ref{thm2}.

\noindent\hspace{2em}{\it Proof of Theorem \ref{thm1}: }
Applying equation (\ref{11}) to the two hypotheses $p_{\theta-\epsilon}$ and $p_{\theta+\epsilon}$, we obtain 
\begin{align*}
\inf_{\theta-\epsilon\le \theta'
\le \theta+\epsilon}
\beta(\vec{T},\theta',\epsilon)\le
\min\{ \beta^+(\vec{T},\theta-\epsilon,\epsilon),
\beta^-(\vec{T},\theta+\epsilon,\epsilon) \}
\le
\sup_{0 \,< s \,< 1} I^s(p_{\theta-\epsilon} \|
p_{\theta+\epsilon}).
\end{align*}
Taking the limit $\epsilon \to 0$, we have
\begin{align}
\lim _{\epsilon \to + 0}
\frac{1}{g(\epsilon)}
\sup_{\vec{T}:WC} 
\inf_{\theta-\epsilon\le \theta'
\le \theta+\epsilon}
\beta(\vec{T},\theta',\epsilon)
\le
2^{\kappa}
\sup_{0 \,< s \,< 1}I^s_{g,\theta}. \Label{hise1}
\end{align}
\endproof
In the regular case, this method was used by Sievers \cite{Sie}.

\noindent\hspace{2em}{\it Proof of Theorem \ref{thm2}:} 
Hoeffding's formula (\ref{Hoe}) yields
\begin{align*}
\inf_{\theta-(1-\eta)\epsilon \le
 \theta'\le \theta+ (1-\eta)\epsilon}
\beta(\vec{T},\theta',(1-\eta)\epsilon)
&\le
\beta^-(\vec{T},\theta+ \frac{1}{2}\epsilon,(1-\eta)\epsilon) \\
&\le
\sup_{0\,< s \,< 1}
\frac{-s 
\beta^+(\vec{T},\theta- \frac{1}{2}\epsilon,\eta\epsilon)
+ I^s(p_{\theta- \frac{1}{2}\epsilon}\|
p_{\theta+ \frac{1}{2}\epsilon})}
{1-s} \\
&\le
\sup_{0\,< s \,< 1}
\frac{-s 
\inf_{\theta-\eta\epsilon \le
 \theta'\le \theta+ \eta\epsilon}
\beta(\vec{T},\theta',\eta\epsilon)
+ I^s(p_{\theta- \frac{1}{2}\epsilon}\|
p_{\theta+ \frac{1}{2}\epsilon})}{1-s}.
\end{align*}
The uniformity of (\ref{uni}) guarantees that
\begin{align}
\alpha_2(\vec{T},\theta)
(1-\eta)^\kappa
&\le
\sup_{0\,< s \,< 1}
\frac{-s \alpha_2(\vec{T},\theta) \eta^\kappa
+ I^s_{g,\theta}}{1-s}.
\Label{1}
\end{align} From (\ref{1}), we have
\begin{align}
\alpha_2(\vec{T},\theta)
\le 
\sup\left\{
\alpha\left|
\alpha ( \eta^\kappa, (1-\eta)^\kappa)
\in
\left\{
(x,y)\left|
y \le \max_{0 \le s\le 1}
\frac{-sx +I^s_g}
{1-s}, x, y \ge 0
\right.
\right\}, \quad
0 \le \forall \eta \le 1
\right.
\right\}. \Label{cut3}
\end{align}
We define the set ${\cal C}_1$ and $\alpha_0$ as
\begin{align}
{\cal C}_1&:=
\left\{
(x,y)\left| y \ge \sup_{0\,< t \,< 1}
\frac{-tx +I^t_{g,\theta}}{1-t}
\right.\right\} \\
\alpha_0&:= \sup\{\alpha|
(\alpha \eta^{\kappa},
\alpha(1-\eta)^{\kappa})
\notin \inner({\cal C}_1),
\quad 0\le \forall \eta \le 1\}. \Label{28-1}
\end{align}
Note that the function $s \mapsto I^s_{g,\theta}$ is concave. We define the convex function $g(x)$ and another set ${\cal C}_2$ as
\begin{align}
g(x)&:=\left\{
\begin{array}{cc}
\alpha_0 \left( 1- \left(\frac{x}{\alpha_0}\right)
^{\frac{1}{\kappa} }\right)^{\kappa}&
0\le x \le \alpha_0 \\
0 & x \,> \alpha_0 
\end{array}
\right. \Label{28-2}\\
{\cal C}_2 &:= \{ (x,y) | y \ge g(x)\}.\Label{28-3}
\end{align}
Since inequality (\ref{cut3}) guarantees 
\begin{align*}
(\alpha_2(\vec{T},\theta)\eta^{\kappa},
\alpha_2(\vec{T},\theta)(1-\eta)^{\kappa})
\notin \inner({\cal C}_1), 
\quad 0\le \forall \eta \le 1,
\end{align*}
the inequality 
\begin{align}
\alpha_2(\vec{T},\theta) \le \alpha_0
\end{align}
holds. Relations (\ref{28-1}), (\ref{28-2}), and (\ref{28-3})
guarantee the 
relation
\begin{align}
{\cal C}_1 \subset {\cal C}_2 \Label{28-4}.
\end{align}
Applying Lemma \ref{concave}, we have
\begin{align}
& \inf_{(x,y)\in {\cal C}_1}
(sx +(1-s)y)
=
\inf_{x \ge 0}
\left(sx +(1-s)\sup_{0 \,< t \,< 1}
\frac{-tx +I^t_{g,\theta}}{1-t}\right) \nonumber\\
=&
\inf_{x \,>0}\sup_{0 \,< t \,< 1}
\left( \frac{(s-t)x +(1-s)I^t_{g,\theta}}{1-t}\right)
= I^s_{g,\theta}. \Label{4-19}
\end{align}

In the following, we divide our situation into three cases $\kappa \,> 1, \kappa =1 , 1 \,> \kappa \,> 0$. When $\kappa \,> 1$, relation (\ref{28-4}) guarantees that 
\begin{align*}
&\alpha_0 s(1-s)\left( s^{\frac{1}{\kappa-1}}
+ (1-s)^{\frac{1}{\kappa-1}}\right)
^{1-\kappa}
=
\alpha_0 \min_{0\,< \eta \,< 1}
\left(s\eta^{\kappa}+(1-s)(1-\eta)^{\kappa}\right)\\
=& \inf_{(x,y) \in {\cal C}_2} \left(sx + (1-s)y \right)
\le  \inf_{(x,y) \in {\cal C}_1} (sx + (1-s)y )
= I^s_{g,\theta}.
\end{align*}
Therefore, 
\begin{align*}
\alpha_2(\vec{T},\theta)
\le
\alpha_0
\le
\frac{ \left( s^{\frac{1}{\kappa-1}}
+ (1-s)^{\frac{1}{\kappa-1}}\right)
^{\kappa-1}}{s(1-s)}
I^s_{g,\theta},
\end{align*}
which implies (\ref{9}).

When $\kappa =1$, similarly, we can easily prove 
\begin{align*}
\frac{1}{2}
\alpha_2(\vec{T},\theta)
=\alpha_2(\vec{T},\theta)
\min_{0\,< \eta\,<1} \left(\frac{1}{2}\eta
+\frac{1}{2}(1-\eta)\right)
\le \inf_{(x,y)\in {\cal C}_1} 
\left(\frac{1}{2} x + \frac{1}{2} y \right)
= I^{\frac{1}{2}}_{g,\theta}.
\end{align*}

Finally, we consider the case where $\kappa \,< 1$. Since function $g$ is concave on $(0,\alpha_0)$, there exists $\eta_0 \in (0,1)$ such that
\begin{align}
\alpha_0 (\eta_0^{\kappa} ,(1-\eta_0)^{\kappa})
\in {\cal C}_1 \cap \overline{{\cal C}_2^c}.
\end{align}
Since ${\cal C}_1$ and ${\cal C}_2^c$ are convex, there exists a real number $s_0 \in (0,1)$ such that 
\begin{align*}
\sup_{(x,y) \in  {\cal C}_2^c}
\left(s_0 x + (1-s_0)y \right)
=
\alpha_0 ( s_0 \eta_0^{\kappa} 
+ (1-s_0)(1-\eta_0)^{\kappa})
=\inf_{(x,y) \in  {\cal C}_1}
\left(s_0 x + (1-s_0)y\right).
\end{align*}
In general, for any $s \in (0,1)$, using (\ref{4-19}), we obtain 
\begin{align*}
& \alpha_0 s(1-s) \left(s^{\frac{1}{\kappa-1}}
+(1-s)^{\frac{1}{\kappa-1}}\right)
^{1-\kappa} 
= \alpha_0 \sup_{0\,< \eta \,< 1}
\left(s \eta^{\kappa} +(1-s) (1-\eta)^{\kappa} \right)
=\sup_{(x,y) \in  {\cal C}_2^c}
\left(s x + (1-s)y \right)\\
\ge & \alpha_0 \left(
s \eta_0^{\kappa} +(1-s) (1-\eta_0)^{\kappa}\right)
\ge 
\inf_{(x,y) \in  {\cal C}_1}
\left(s x + (1-s)y\right) =
I^{s}_{g,\theta},
\end{align*}
which lead to (\ref{9}).
\endproof

\noindent\hspace{2em}{\it Proof of Theorem \ref{thm3}:}
It is trivial in the case of $\kappa=1$.
When $\kappa \,>1$, it is also trivial because 
\begin{align*}
\inf_{0 \,< s \,< 1}
\frac{I^s_{g,\theta}}{s(1-s)}
\left( s^{\frac{1}{\kappa-1}} 
+ (1-s)^{\frac{1}{\kappa-1}} 
\right)^{\kappa-1}
\le
2^\kappa I_{g,\theta}^{\frac{1}{2}}
\le
2^\kappa \sup_{0\le s \le 1}
I_{g,\theta}^{s}.
\end{align*}
Next, we consider the case $\kappa \,< 1$. The inequality
\begin{align}
\sup_{0\,< s \,< 1}
\frac{I^s_{g,\theta}}{s(1-s)}
\left( s^{\frac{1}{\kappa-1}} 
+ (1-s)^{\frac{1}{\kappa-1}} 
\right)^{\kappa-1}
\le
2^\kappa \sup_{0\le s \le 1}
I_{g,\theta}^{s} \Label{19}
\end{align}
follows from the two inequalities
\begin{align}
I^s_{g,\theta}
&\le \sup_{0 \,< s \,< 1}
I_{g,\theta}^{s} \Label{17}\\
\frac{1}{s(1-s)}
\left( s^{\frac{1}{\kappa-1}} 
+ (1-s)^{\frac{1}{\kappa-1}} 
\right)^{\kappa-1}
&\le
2^\kappa \Label{18.9}.
\end{align}
We thus obtain (\ref{14}). 
In the following, we prove that the equality of (\ref{14}) implies (\ref{16.3}) and (\ref{16.3}) implies (\ref{15}) and the equality of (\ref{14}) in the case where $\kappa \,< 1$. If we assume that the equality of (\ref{19}) holds, the equalities of (\ref{17}) and (\ref{18.9}) hold at the same $s$. The equality of (\ref{18.9}) holds 
if and only if $s=\frac{1}{2}$. Therefore, 
\begin{align*}
I^{\frac{1}{2}}_{g,\theta}
= \sup_{0 \,< s \,< 1}
I_{g,\theta}^{s},
\end{align*}
which is equivalent to (\ref{16.3}). If we assume that (\ref{16.3}) holds, inequality (\ref{18.9}) guarantees that
\begin{align}
\overline{\alpha}_2(\theta) 
\le 2^{\kappa} I^{\frac{1}{2}}_{g,\theta}.
\end{align}
Substituting $\frac{1}{2}$ into $s$ at the 
left hand side (LHS) in the definition of $\overline{\alpha}_2(\theta)$, we obtain
\begin{align*}
\overline{\alpha}_2(\theta) \ge
2^{\kappa} I^{\frac{1}{2}}_{g,\theta}.
\end{align*}
Thus, equation (\ref{15}) holds. Combining (\ref{16.3}) and (\ref{15}), we obtain the equality of (\ref{14}).
\endproof
\section{First exponential rates of useful estimators}
\Label{s5}
In the following, we discuss the first exponential rate $\beta(\vec{T}, \theta,\epsilon)$ of a useful estimator $\vec{T}$ for a location shift family $\{ f(x - \theta)| \theta \in \real\}$, where $f$ is a probability density function on $\real$.

When the function $d_{\epsilon}(x):=f(x + \epsilon) /f(x -\epsilon)$ is monotonically decreasing w.r.t. $x$ when both $f(x + \epsilon)$ and $f(x -\epsilon)$ are not zero, we can define the likelihood ratio estimator $\vec{\tilde{\theta}_{\epsilon}}:
=\{\tilde{\theta}_{n,\epsilon}(x_1 , \ldots , x_n)\}$,
which depends on the constant $\epsilon \,> 0$, as shown by
\begin{align}
\tilde{\theta}_{n,\epsilon}:=
\frac{1}{2} \left(\sup\{ z | k(z) \,< 0 \}+
\inf \{ z | k(z) \,> 0 \} \right), \Label{725.2}
\end{align}
where the monotonically decreasing function $k(z)$ is defined by
\begin{align}
k(z):= \frac{1}{n} \sum_{i=1}^n \left(
\log f( x_i - z +\epsilon ) - 
\log f( x_i - z -\epsilon )\right)
\Label{18.1}.
\end{align}
Note that when $\log f(x)$ is concave, the above condition is satisfied. This definition is well defined although the monotonically decreasing function $k(z)$ is not continuous.

If the support of $f$ is $(a,b)$, we need to modify the definition as follows. In this case, we modify the estimator $\tilde{\theta}_{n,\epsilon}$ by using the two estimators $\overline{\theta}_n:= \max\{ x_1, \ldots, x_n \} - b$ and $\underline{\theta}_n:= \min\{ x_1, \ldots, x_n \} - a$. When $\underline{\theta}_n- \overline{\theta}_n \,> 2 \epsilon$, the estimated value is defined by (\ref{725.2}) in the interval $(\underline{\theta}_n- \epsilon ,\overline{\theta}_n + \epsilon)$. When $\underline{\theta}_n- \overline{\theta}_n \le 2 \epsilon$, we define $\tilde{\theta}_{n,\epsilon}:= \frac{1}{2}\left(\underline{\theta}_n+ \overline{\theta}_n\right)$. Moreover, when the support of $f$ is a half line $(0, \infty)$, the estimated value is defined by (\ref{725.2}) in the half line $(\underline{\theta}_n- \epsilon , \infty)$.

We have the following lemma. 
(A regular version of this lemma is discussed by Huber \cite{Hu}, Sievers \cite{Sie}, and Fu \cite{Fu85}.)
\begin{lem}\Label{thm7.2}
When $\log f(x)$ is concave, the equation 
\begin{align}
\min\{ 
\beta^-( \vec{\tilde{\theta}_{\epsilon}}, 
\theta , \epsilon ),
\beta^+( \vec{\tilde{\theta}_{\epsilon}}, 
\theta , \epsilon )\}
= \sup_{0\,< s \,< 1}
I^s(f_{\theta-\epsilon }\| f_{\theta+\epsilon}) . 
\Label{725.3}
\end{align}
holds, where $f_{\theta}(x):= f(x- \theta)$.
\end{lem}
Therefore, the equality of inequality (\ref{7}) holds in this case.

\begin{proof}
Note that 
$\beta^+( \vec{\tilde{\theta}_{\epsilon}}, 
\theta , \epsilon )
=\beta^+( \vec{\tilde{\theta}_{\epsilon}}, 
\theta - \epsilon , 
\epsilon)$
and 
$\beta^-( \vec{\tilde{\theta}_{\epsilon}},
 \theta , \epsilon )
= \beta^-( \vec{\tilde{\theta}_{\epsilon}}, 
\theta + \epsilon , 
\epsilon )$
because of the shift-invariance. From the concavity, the condition $\tilde{\theta}_{\epsilon,n} \le \theta$ is equivalent to the condition $\sup\{ z | k(z) \,< 0\} \le \theta$, which implies that $k(\theta)\ge 0$.
Thus, we have 
\begin{align*}
\sum_i
\log f(x_i - \theta +\epsilon)
- \log f(x_i - \theta -\epsilon)
\ge 0.
\end{align*}
Conversely,
the condition \begin{align*}
\sum_i
\log f(x_i - \theta +\epsilon)
- \log f(x_i - \theta -\epsilon)
\,> 0
\end{align*}
implies that
$k(\theta) \,> 0$.
Thus, we have
$\sup\{ z | k(z) \,< 0\} \le \theta$,
which is equivalent to the condition $\tilde{\theta}_{\epsilon,n} \le \theta$. Therefore, we have the relations 
\begin{align*}
\{ f_{\theta-\epsilon}^n (\vec{x}_n) \,>
f_{\theta+\epsilon}^n (\vec{x}_n) \}
\subset 
\{ \tilde{\theta}_{n,\epsilon} \le \theta\}
\subset
\{ f_{\theta-\epsilon}^n (\vec{x}_n) \ge
f_{\theta+\epsilon}^n (\vec{x}_n) \}. 
\end{align*}
Similarly, we can prove
\begin{align*}
\{ f_{\theta+\epsilon}^n (\vec{x}_n) \,>
f_{\theta-\epsilon}^n (\vec{x}_n) \}
\subset 
\{ \tilde{\theta}_{n,\epsilon} \ge \theta\}
\subset
\{ f_{\theta+\epsilon}^n (\vec{x}_n) \ge
f_{\theta-\epsilon}^n (\vec{x}_n) \} .
\end{align*}
Applying (\ref{11}), we can prove
\begin{align*}
\min\{ 
\beta^-( \vec{\tilde{\theta}_{\epsilon}}, 
\theta + \epsilon, 
\epsilon ),
\beta^+( \vec{\tilde{\theta}_{\epsilon}}, 
\theta - \epsilon,
 \epsilon )\}
= \sup_{0\,< s \,< 1}
I^s(f_{\theta-\epsilon }\| f_{\theta+\epsilon}),
\end{align*}
which implies equation (\ref{725.3}).
\end{proof}

\begin{cor}
When the function $\log f(x)$ is concave, the equation 
\begin{align*}
\alpha_1(\theta)= \overline{\alpha}_1(\theta)
\end{align*}
holds.
\end{cor}
\begin{lem}
When $f(x)$ is monotonically decreasing, the estimator $\vec{\underline{\theta}_\epsilon} := \{ \underline{\theta}_{\epsilon,n} := \underline{\theta}_n - \epsilon\}$ satisfies the relations
\begin{align}
\beta^+(\vec{\underline{\theta}_\epsilon}, 
\theta, \epsilon)
&=\sup_{0\,< s \,< 1} I^s (f_{\theta-\epsilon}
\|f_{\theta+\epsilon}) \nonumber \\
\beta^-(\vec{\underline{\theta}_\epsilon}, 
\theta, \epsilon)
&=\infty \Label{3-8-1}.
\end{align}
Thus, in this case, the equation 
\begin{align*}
\alpha_1(\theta)= \overline{\alpha}_1(\theta)
\end{align*}
holds.
\end{lem}
\begin{proof}
Since $\underline{\theta}_n \,> \theta$, we have $\underline{\theta}_{\epsilon,n} \,> \theta- \epsilon$, which implies (\ref{3-8-1}). If $\underline{\theta}_{\epsilon,n} \ge \theta +\epsilon$, we have $\underline{\theta}_{n} \ge \theta + 2 \epsilon$.
Thus, $f(x_i -(\theta +2\epsilon)) 
\ge f(x_i -\theta )$ for any $i= 1, \ldots, n$.
Therefore, 
\begin{align*}
f^n_{\theta+2\epsilon} (\vec{x}_n)
\ge f^n_{\theta} (\vec{x}_n).
\end{align*}
Conversely, if $f^n_{\theta+2\epsilon} (\vec{x}_n) \ge f^n_{\theta} (\vec{x}_n)$,
we have $\underline{\theta}_{n} \ge \theta + 2 \epsilon$. Thus, 
\begin{align*}
f^n_{\theta+2\epsilon} \{ \vec{x}_n|
f^n_{\theta+2\epsilon} (\vec{x}_n)
\,< f^n_{\theta} (\vec{x}_n)\}
= 
f^n_{\theta+2\epsilon} \{ \vec{x}_n|
\underline{\theta}_{n} 
\,< \theta + 2 \epsilon\}=0.
\end{align*}
Since the likelihood test $\{ \vec{x}_n| \underline{\theta}_{n} \ge \theta + 2 \epsilon\}$ achieves the optimal rate (\ref{11}), we have 
\begin{align*}
&\lim -\frac{1}{n}\log
f^n_{\theta} \{ \vec{x}_n|
f^n_{\theta+2\epsilon} (\vec{x}_n)
\ge f^n_{\theta} (\vec{x}_n)\} \\
=&
\lim -\frac{1}{n}\log
\left(
f^n_{\theta} \{ \vec{x}_n|
f^n_{\theta+2\epsilon} (\vec{x}_n)
\ge f^n_{\theta} (\vec{x}_n)\}+
f^n_{\theta+2\epsilon} \{ \vec{x}_n|
f^n_{\theta+2\epsilon} (\vec{x}_n)
\,< f^n_{\theta} (\vec{x}_n)\}\right)\\
=&
\sup_{0\,< s \,<1}
I^s(f_{\theta}\|f_{\theta+2\epsilon}) 
= 
\sup_{0\,< s \,<1}
I^s(f_{\theta- \epsilon}\|f_{\theta+\epsilon}) .
\end{align*}
\end{proof}
\begin{lem}\Label{l10}
When the function $x \mapsto \log f(x)$ is concave, MLE $\vec{\theta_{ML}}:=\{\theta_{ML,n}\}$ satisfies that 
\begin{align}
\beta^+(\vec{\theta_{ML}}, \theta,\epsilon) 
&= \sup_{t \ge 0}
- \log \int_{a+\epsilon}^b \exp \left( -t 
\frac{f'(x-\epsilon)}{f(x-\epsilon)} \right)
f(x) \,d x \Label{33} \\
&= \sup_{t \ge 0}
- \log \int_{a}^{b-\epsilon} \exp \left( -t 
\frac{f'(x)}{f(x)} \right)
f(x+\epsilon) \,d x \Label{34}\\
\beta^-(\vec{\theta_{ML}}, \theta,\epsilon) 
&= \sup_{t \ge 0}
- \log \int_{a}^{b-\epsilon} \exp \left( t 
\frac{f'(x+\epsilon)}{f(x+\epsilon)} \right)
f(x) \,d x \Label{35}\\
&= \sup_{t \ge 0}
- \log \int_{a+\epsilon}^b \exp \left( t 
\frac{f'(x)}{f(x)} \right)
f(x-\epsilon) \,d x \Label{36}.
\end{align}
\end{lem}
This lemma is a special case of Fu's result \cite{Fu73}. 
For the reader's convenience, we give its proof.

\begin{proof}
Equations (\ref{34}) and (\ref{36}) are trivial. We prove (\ref{33}). From the assumption that for any $\vec{x}_n:= (x_1 , \ldots, x_n)$, the function $\theta \mapsto \sum_{i =1}^n \log f(x_i - \theta)$ is concave 
on $(\overline{\theta}(\vec{x}_n), \underline{\theta}(\vec{x}_n))$, the function $\theta \mapsto \sum_{i=1}^n l_\theta(x_i)$ is monotonically decreasing, 
where $l_\theta(x):=-\frac{f'(x-\theta)}{f(x-\theta)}$ on $(\overline{\theta}(\vec{x}_n), \underline{\theta}(\vec{x}_n))$. 
When $\vec{x}_n$ belongs to the support of $f_{\theta'}$, the condition $\theta_{ML,n}(\vec{x}_n) \ge \theta'$ is equivalent to the condition 
\begin{align*}
\frac{1}{n} \sum_{i=1}^n l_{\theta'}(x_i) \ge 0.
\end{align*}
Denoting the conditional probability $f\{ A| x \in B\}$ under the condition $x \in B$, we can evaluate 
\begin{align}
\lim - \frac{1}{n} \log f_\theta^n 
\{ \theta_{ML,n} \ge \theta+\epsilon\}
&=
\lim -\frac{1}{n}\log f_{\theta}^n \{
\theta_{ML,n} \ge \theta+\epsilon |
\vec{x}_n \in (a+\epsilon,b+\epsilon)^n\}
- \frac{1}{n}
\log f^n_{\theta}(a+\epsilon,b+\epsilon)^n 
\nonumber\\
&= \lim -\frac{1}{n}\log f_{\theta,\epsilon}^n 
\left \{
\frac{1}{n} \sum_{i=1}^n l_{\theta+\epsilon}(x_i) 
\ge 0 \right\}
-\log \int_{a+\epsilon}^b f(x) \,d x , \Label{31}
\end{align}
where the probability density function whose support is $(a+\epsilon,b)$ is defined by
\begin{align*}
f_{\theta,\epsilon}(x):=
\frac{f(x)}{\int_{a+\epsilon}^b f(x) \,d x}.
\end{align*}
Chernoff's theorem (Theorem 3.1 in Bahadur \cite{Ba71}) guarantees that
\begin{align} 
\lim -\frac{1}{n}\log f_{\theta,\epsilon}^n 
\{
\frac{1}{n} \sum_{i=1}^n l_{\theta+\epsilon}(x_i) 
\ge 0 \}
&=
\sup_{t \ge 0}
- \log \int_{a+\epsilon}^b \exp \left(
t l_{\theta+\epsilon}(x) \right)
f_{\theta,\epsilon}(x) \,d x \nonumber\\
&=
\sup_{t \ge 0}
- \log \int_{a+\epsilon}^b \exp \left(
-t \frac{f'(x-\theta)}{f(x-\theta)}\right)
f(x) \,d x
+\log \int_{a+\epsilon}^b f(x) \,d x .\Label{32}
\end{align}
Combining (\ref{31}) and (\ref{32}), we obtain (\ref{33}). Similarly, We can prove (\ref{35}).
\end{proof}
\begin{lem}
When $f(x)$ is monotonically decreasing, the MLE $\theta_{ML,n}$ equals the estimator $\underline{\theta}_n$.
\end{lem}
\begin{proof}
For any data $\vec{x}_n:= (x_1, \ldots, x_n)$, if $\theta \,< \underline{\theta}_n(\vec{x}_n))$, $f(x_1- \theta) \ldots f(x_n- \theta)=0$. Conversely, if $\theta \,> \underline{\theta}_n(\vec{x}_n))$, the obtained $f(x_i- \theta) \le f(x_i - \underline{\theta}_n(\vec{x}_n))$. Thus, $\underline{\theta}_n$ is the MLE. 
\end{proof}
\begin{lem}\Label{l34.3}
Let $f$ be a density function whose support is the interval $(a,b)$.
The estimators $\vec{\underline{\theta}}$ and $\vec{\overline{\theta}}$ satisfy 
\begin{align}
\beta^+(\vec{\underline{\theta}},\theta,\epsilon)
&=
- \log \left(\int_{a}^{b-\epsilon}f(x) \,d x\right)
,\quad
\beta^-(\vec{\underline{\theta}},\theta,\epsilon)
=\infty \Label{3-8-14} \\
\beta^+(\vec{\overline{\theta}},\theta,\epsilon)
&=\infty ,\quad
\beta^-(\vec{\overline{\theta}},\theta,\epsilon)
=
- \log \left(\int_{a+\epsilon}^{b}f(x) \,d x \right)
\Label{3-8-15}
\end{align}
We can use the convex combination $\vec{\check{\theta}(\lambda)} := \{ \check{\theta}(\lambda)_n:= \lambda \underline{\theta}_n + (1-\lambda) \overline{\theta}_n\}$ with the ratio $\lambda : 1- \lambda$, where $0 \,< \lambda \,< 1$. It satisfies that 
\begin{align}
\beta^+(\vec{\check{\theta}(\lambda)},\theta,\epsilon)
&=
- \log \left(\int_{a}^{b-\frac{\epsilon}{1-\lambda}}f(x) \,d x\right)
\Label{34.1}\\
\beta^-(\vec{\check{\theta}(\lambda)},\theta,\epsilon)
&=
- \log \left(\int_{a+\frac{\epsilon}{\lambda}}^{b}f(x) \,d x \right)
\Label{34.2}.
\end{align}
\end{lem}
\begin{proof}
Define $\overline{\omega}_n:=\max \{ \omega_1 , \ldots , \omega_n \}, \underline{\omega}_n:=\min \{ \omega_1 , \ldots , \omega_n \}$. Since the estimators $\vec{\underline{\theta}}, \vec{\overline{\theta}}$, and $\vec{\check{\theta}(\lambda)}$ are covariant for location shift, we may discuss only the case that $\theta=0$. From the relation $\underline{\theta}_n \,> \theta \,> \theta - \epsilon$, we obtain the second equation of (\ref{3-8-14}). Its joint probability density function $f_n( \overline{\omega}_n,\underline{\omega}_n)$ is given by
\begin{align*}
f_n( \overline{\omega}_n,\underline{\omega}_n)
:=
\begin{cases}
n(n-1)\left( \int_{\underline{\omega}_n}^{\overline{\omega}_n}
f(x) \,d x \right)^{n-2}f(\underline{\omega}_n)f(\overline{\omega}_n)
& \overline{\omega}_n \ge \underline{\omega}_n \\
0& \overline{\omega}_n \,< \underline{\omega}_n .
\end{cases}
\end{align*}
Defining
\begin{align*}
g(\underline{\omega}_n,\overline{\omega}_n):=
\begin{cases}
\int_{\underline{\omega}_n}^{\overline{\omega}_n}
f(x) \,d x & \overline{\omega}_n \ge \underline{\omega}_n \\
0& \overline{\omega}_n \,< \underline{\omega}_n ,
\end{cases}
\end{align*}
we have 
\begin{align}
p^n_{\theta}(\check{\theta}(\lambda)_n \le \theta- \epsilon)
= \int_{ \check{\theta}(\lambda)
(\underline{\omega}_n, \overline{\omega}_n) \le \epsilon}
n(n-1) g( \underline{\omega}_n,\overline{\omega}_n)^n
f(\underline{\omega}_n)f(\overline{\omega}_n)
\,d \underline{\omega}_n
\,d \overline{\omega}_n \Label{34.6}.
\end{align} From the continuity of 
$f(\underline{\omega}_n), f(\overline{\omega}_n)$ and
$g( \underline{\omega}_n,\overline{\omega}_n)$,
the equations 
\begin{align*}
\lim_{n \to \infty}
\frac{1}{n}\log p^n_{\theta}(\check{\theta}(\lambda)_n
\le \theta- \epsilon)
&= 
\lim_{n \to \infty}\frac{1}{n}
\sup_{\check{\theta}(\lambda)(\underline{\omega}_n, \overline{\omega}_n)  \le \epsilon}
\log g(\underline{\omega}_n,\overline{\omega}_n) \\
&= - \log \left(\int_{a}^{b-\frac{\epsilon}{1-\lambda}}f(x) \,d x\right)
\end{align*}
hold.
This implies the first equation of (\ref{3-8-14}) and (\ref{34.1}).
In addition, we can similarly show the same for (\ref{3-8-15}) and (\ref{34.2}).
\end{proof}

\section{Two bounds in location shift family}
\Label{s6}
We discuss a location shift family generated by a probability density function (pdf) whose support is an interval $(a,b)$. Moreover, we assume that the pdf $f$ is $C^1$ continuous and satisfies that
\begin{align*}
f(x) \cong A_1 (x-a)^{\kappa_1-1}, & \quad x \to a +0\\
f(x) \cong A_2 (b- x)^{\kappa_2-1},& \quad x \to b -0,
\end{align*}
where $\kappa_1,\kappa_2 \,> 0$, as for the uniform and beta distributions.

When its support is a half line $(0,\infty)$ as for the gamma distribution and Weibull distribution, our situation results in the above case where $A_2=0$ if $f$ is $C^3$ continuous and $\lim_{x \to \infty}
|\frac{\,d }{\,d x} \log f(x)|\,< \infty$.
Also, when $\kappa_1 \,> \kappa_2$, our situation results in the above case where $A_2=0$.

In the following, in this setting, we calculate two upper bounds $\overline{\alpha}_1(\theta)$ and $\overline{\alpha}_2(\theta)$, and derive a necessary and sufficient condition for coincidence of the two upper bounds. In addition, in some cases, we calculate the two bounds $\alpha_1(\theta)$ and $\alpha_2(\theta)$.
\subsection{Semi-regular case}\Label{e1}
As was proven in \cite{Ak}, when $\kappa_1,\kappa_2 \,> 2$, the relation
\begin{align}
\lim_{\epsilon \to +0}
\frac{I^s( f_{\theta} \| f_{\theta+ \epsilon})}
{\epsilon^2}
=
\frac{s(1-s)}{2} J_f \Label{46.7}
\end{align}
holds, where this convergence is uniform for $s$ and $J_f$ is defined by
\begin{align*}
J_f:= \int_a^b 
\left(\frac{\,d f(x)}{\,d x}\right)^2 f(x)^{-1}
\,d x.
\end{align*}
\begin{prop}
When $g(x)=x^2$, we obtain $\kappa=2$ and the relation
\begin{align}
\overline{\alpha}_1(\theta)=\alpha_1(\theta)=
\overline{\alpha}_2(\theta)
=\alpha_2(\theta)=
\frac{1}{2} J_\theta. \Label{39}
\end{align}
\end{prop}
\begin{proof}
Using (\ref{46.7}), we have
\begin{align*}
\overline{\alpha}_1(\theta) &= 
4 \max_{0\le s \le 1}
I^s_{g,\theta}
=\frac{1}{2} J_\theta \\
\overline{\alpha}_2(\theta) &= 
\min_{0\le s \le 1}
\frac{I^s_{g,\theta}}{s(1-s)}
=\frac{1}{2} J_\theta .
\end{align*}
When the function $\log f(x)$ is concave, using Lemma \ref{l10} we can evaluate
\begin{align*}
\beta^+(\vec{\theta_{ML}}, \theta,\epsilon) 
& \ge 
- \log \int_{a}^{b-\epsilon} \exp \left( -\epsilon 
\frac{f'(x)}{f(x)} \right)
f(x+\epsilon) \,d x \\
& \cong - \log \left( \int_a^{b-\epsilon}
\left( 1- \epsilon \frac{f'(x)}{f(x)}
+ \frac{\epsilon^2}{2}
\left(\frac{f'(x)}{f(x)}\right)^2
\right)
\left( f(x) + \epsilon f'(x)
+ \frac{\epsilon^2}{2}f''(x)\right)
\,d x + o( \epsilon^2) \right) \\
& \cong - \log \left( \int_a^{b-\epsilon}
f(x) - \frac{\epsilon^2}{2} 
\left(\frac{f'(x)^2}{f(x)} + f''(x)\right)
\,d x + o( \epsilon^2) \right) \\
& \cong - \log \left( 
1- 
\int_{b-\epsilon}^b f(x) \,d x 
- \frac{\epsilon^2}{2} \int_a^{b-\epsilon}
\frac{f'(x)^2}{f(x)}\,d x + o( \epsilon^2) \right) \\
& \cong - \log \left( 
1 
- \frac{\epsilon^2}{2} 
\left(\int_a^{b-\epsilon}
\frac{f'(x)^2}{f(x)}\,d x 
-f'(b-\epsilon)\right)
+ o( \epsilon^2) \right) .
\end{align*}
Thus, we obtain 
\begin{align*}
\lim_{\epsilon\to 0} \frac{1}{\epsilon^2}
\beta^+(\vec{\theta_{ML}}, \theta,\epsilon)
\ge \frac{1}{2} \lim_{\epsilon\to 0}
\left(\int_a^{b-\epsilon}
\frac{f'(x)^2}{f(x)}\,d x -f'(b-\epsilon)\right)
= \frac{1}{2}J_f.
\end{align*}
Similarly, we can prove that
\begin{align*}
\lim_{\epsilon\to 0} \frac{1}{\epsilon^2}
\beta^-(\vec{\theta_{ML}}, \theta,\epsilon)\ge \frac{1}{2}J_f.
\end{align*}
Since $\alpha_2(\vec{\theta_{ML}}, \theta)\le\frac{1}{2}J_f$, we have
\begin{align*}
\alpha_2(\vec{\theta_{ML}}, \theta)
=  \frac{1}{2}J_f,
\end{align*}
which implies (\ref{39}).
\end{proof}
\subsection{The case that $\kappa_1=\kappa_2=1$}
\Label{e2}
As was proven elsewhere \cite{Haya2} \cite{Ak},
when $\kappa_1=\kappa_2=1$, the equation
\begin{align*}
\lim_{\epsilon \to +0}
\frac{I^s( f_{\theta} \| f_{\theta+ \epsilon}) }
{\epsilon}
= A_1s + A_2 (1-s)
\end{align*}
holds, where this convergence is uniform for $s \in (0,1)$.
Letting $g(x)=|x|$, we have $\kappa=1$.
\begin{prop}
The relations
\begin{align}
\alpha_1(\theta) 
&= \overline{\alpha}_1(\theta) = 
2 \sup_{0\,< s \,< 1}
I^s_{g,\theta}
=2 \max \{A_1, A_2\} \Label{3-8-2}\\
\alpha_2(\theta) &= 
\overline{\alpha}_2(\theta) = 
2 I^{\frac{1}{2}}_{g,\theta}=
A_1+A_2 \Label{3-8-3}
\end{align}
hold.
Therefore, $\alpha_1(\theta) = \alpha_2(\theta) $ 
if and only if $A_1=A_2$.
\end{prop}
\begin{proof}
The third equations of (\ref{3-8-2}) and (\ref{3-8-3}) follow from the formula $I^s_{g,\theta}= A_1s + A_2(1-s)$. In the following, we prove the first equations of (\ref{3-8-2}) and (\ref{3-8-3}).
By using Lemma \ref{l34.3}, the equations 
\begin{align*}
\lim_{\epsilon \to +0}
\frac{1}{\epsilon}
\beta^+(\vec{\check{\theta}(\lambda)}, 
\theta, \epsilon )
=
\frac{A_1}{\lambda} , \quad
\lim_{\epsilon \to +0}
\frac{1}{\epsilon}
\beta^-(\vec{\check{\theta}(\lambda)}, 
\theta, \epsilon )
=
\frac{A_2}{1-\lambda} 
\end{align*}
hold.
Letting $\lambda_0 := A_1/(A_1+A_2)$,
we have 
\begin{align*}
\alpha_2(\vec{\check{\theta}(\lambda_0)}, \theta)
= A_1+A_2,
\end{align*}
which implies the first equation of (\ref{3-8-3}).

Next, we prove the first equation of (\ref{3-8-2}) in the case where $A_1 \ge A_2$. The estimator $\vec{\underline{\theta}_{\epsilon}} :=\{ 
\underline{\theta}_{\epsilon,n}\}$ satisfies 
\begin{align*}
\beta^+(\vec{\underline{\theta}_{\epsilon}},\theta,\epsilon)
= - \log \int_{a+ 2 \epsilon}^b f(x) \,d x, \quad
\beta^-(\vec{\underline{\theta}_{\epsilon}},
\theta,\epsilon)
= \infty.
\end{align*}
Therefore, 
\begin{align*}
\lim_{\epsilon \to +0}
\frac{1}{\epsilon} 
\beta( \vec{\underline{\theta}_{\epsilon}},
\theta,\epsilon)
= 2 A_1,
\end{align*}
which implies the first equation of (\ref{3-8-2}). When $A_2 \ge A_1$, we can similarly prove it.
\end{proof}
\subsection{The case that $\kappa_1=\kappa_2=2$}
\Label{e3}
As was proven by \cite{Haya2}, when $\kappa_1=\kappa_2=2$ the equation 
\begin{align}
\lim_{\epsilon \to +0}
\frac{
I^s( f_{\theta} \| f_{\theta+\epsilon}) }
{-\epsilon^2 \log |\epsilon|}
=\frac{(A_1+A_2) s(1-s)}{2} \Label{15.04}
\end{align} 
holds, where this convergence is uniform for $s \in (0,1)$. Letting $g(x)= - x^2 \log x$, we have $\kappa=2$. Using (\ref{15.04}), we have
\begin{align}
\overline{\alpha}_1(\theta) &= 
4 \sup_{0\,< s \,< 1}
I^s_{g,\theta}
=\frac{A_1+A_2}{2} \Label{3-8-19}\\
\overline{\alpha}_2(\theta) &= 
\inf_{0\,< s \,< 1}
\frac{I^s_{g,\theta}}{s(1-s)}
=\frac{A_1+A_2}{2} . \nonumber
\end{align}
\begin{prop}
When the function $x \mapsto \log f(x)$ is concave, we have
\begin{align}
\alpha_1(\theta) = 
\overline{\alpha}_1(\theta) = 
\alpha_2(\theta) = 
\overline{\alpha}_2(\theta) = 
\frac{A_1+A_2}{2}.\Label{3-8-17}
\end{align}
\end{prop}
\begin{proof}
When $a+\delta \,< x \,< b - \delta$,
we can approximate that
\begin{align}
\exp\left ( - \epsilon 
\frac{f'(x)}{f(x)}\right)
f(x+\epsilon) 
\cong f(x) + \frac{\epsilon^2}{2}
\left( f''(x) - \frac{(f'(x))^2}{f(x)}
\right).
\end{align}
Therefore, from Lemma \ref{l10}, we can evaluate
\begin{align}
& \beta^+(\vec{\theta_{ML}}, \theta,\epsilon) \nonumber \\
&\ge -\log \left(
\int_{a}^{b-\epsilon}
\exp\left ( - \epsilon 
\frac{f'(x)}{f(x)}\right)
f(x+\epsilon) \,d x \right)\nonumber\\
&=-\log \left(
\int_{a+\delta}^{b-\delta}
\exp\left ( - \epsilon 
\frac{f'(x)}{f(x)}\right)
f(x+\epsilon) \,d x +
\int_0^\delta \exp 
\left( -\epsilon
\frac{A_1}{A_1 x}\right)
A_1(x+\epsilon) \,d x +
\int_\epsilon^\delta \exp \left( +\epsilon
\frac{A_2}{A_2 x}\right)A_2(x-\epsilon) \,d x \right)
\nonumber\\
&\cong -\log \Biggl(
\int_{a+\delta}^{b-\delta}
f(x) \,d x
+ \int_{a}^{a+\delta} 
f(x) \,d x + \int_{b- \delta}^{b} 
f(x) \,d x \nonumber\\
& \quad + \quad A_1 \int_0^\delta \left(\exp \left( -\epsilon
\frac{1}{x}\right)(x+\epsilon) - x \right)\,d x +
A_2 \int_\epsilon^\delta \left(\exp \left( +\epsilon
\frac{1}{x}\right)
(x-\epsilon) - x \right) \,d x 
+ o(- \epsilon^2\log \epsilon)\Biggr) \nonumber\\
&\cong 
-\log \left( 1 +\frac{A_1+A_2}{2}\epsilon^2 \log \epsilon +
o(-\epsilon^2 \log \epsilon) \right) \Label{28-10}\\
&\cong -\frac{A_1+A_2}{2}\epsilon^2 \log \epsilon +
o(-\epsilon^2 \log \epsilon),\nonumber
\end{align}
where the relation (\ref{28-10}) follows from Lemma \ref{l11}.
Thus, we have
\begin{align}
\alpha_2(\vec{\theta_{ML}}, \theta) = \frac{A_1+A_2}{2} .
\Label{3-8-18}
\end{align}
Using (\ref{3-8-19}) and (\ref{3-8-18}), we obtain (\ref{3-8-17}).
\end{proof}
\subsection{The case that $1 \,< \kappa_1=\kappa_2 \,< 2$}\Label{e4}
As was proven by Hayashi \cite{Haya2}, when $1 \,< \kappa_1=\kappa_2 \,< 2$, the equation 
\begin{align*}
&\lim_{\epsilon \to +0}
\frac{1}{\epsilon^{\kappa_1}}
I^s( f_{\theta} \| f_{\theta+\epsilon}) \\
=&
\frac{A_1 s(1-s(\kappa_1-1))B(s+\kappa_1(1-s),2-\kappa_1)}{\kappa_1}
+\frac{A_2(1-s)(1-(1-s)(\kappa_1-1))B(1-s+\kappa_1 s,2-\kappa_1)
}{\kappa_1} 
\end{align*}
holds, where this convergence is uniform for $s \in (0,1)$, and $B(x,y)$ is a beta function. Letting $g(x)= |x|^{\kappa_1}$, we have $\kappa=\kappa_1$.
\begin{align*}
&\overline{\alpha}_1(\theta)= 
2^\kappa 
\sup_{0\,< s \,< 1}
I^s_{g,\theta} \\
& =
2^\kappa 
\max_{0\le s \le 1}
\frac{
A_1 s(1-s(\kappa-1))B(s+\kappa(1-s),2-\kappa)
+
A_2(1-s)(1-(1-s)(\kappa-1))
B(1-s+\kappa s,2-\kappa)}
{\kappa} \\
&\overline{\alpha}_2(\theta)= 
\inf_{0\,< s \,< 1}
\frac{I^s_{g,\theta}}{s(1-s)}
\left( s^{\frac{1}{\kappa-1}} 
+ (1-s)^{\frac{1}{\kappa-1}} \right)^{\kappa-1}\\
&=
\inf_{0\,< s \,< 1}
\Biggl[ \frac{
A_1 s(1-s(\kappa-1))B(s+\kappa(1-s),2-\kappa)
+
A_2(1-s)(1-(1-s)(\kappa-1))
B(1-s+\kappa s,2-\kappa)}
{\kappa s(1-s)}\\
&\quad \left( s^{\frac{1}{\kappa-1}} 
+ (1-s)^{\frac{1}{\kappa-1}} \right)^{\kappa-1} 
\Biggr]\\
& \le 
\frac{1}{\kappa}(A_1+A_2) 
2^{\kappa-2}(1-\kappa)
B\left(\frac{1 + \kappa}{2},2 -\kappa\right) .
\end{align*} From Lemma \ref{l34.3},
we have
\begin{align}
\lim_{\epsilon \to +0}
\frac{1}{\epsilon^{\kappa_1}}
\beta^+(\vec{\check{\theta}(\lambda)}, \theta, \epsilon )
&=
A_1\frac{1}{\kappa_1 \lambda^{\kappa_1}} \Label{46.1} \\
\lim_{\epsilon \to +0}
\frac{1}{\epsilon^{\kappa_2}}
\beta^-(\vec{\check{\theta}(\lambda)}, \theta, - \epsilon )
&=
A_2\frac{1}{\kappa_2 (1-\lambda)^{\kappa_2}} \Label{46.2}.
\end{align}
Thus, the estimator 
$\check{\theta}(\lambda)$ achieves the optimal order. However, it does not achieve the optimal coefficient $\overline{\alpha}_2(\theta) $.

\begin{prop}
If, and only if, $A_1= A_2$, the equality 
\begin{align*}
\overline{\alpha}_1(\theta)
=\overline{\alpha}_2(\theta) 
\end{align*}
holds.
In this case, 
\begin{align}
\overline{\alpha}_1(\theta)
=\overline{\alpha}_2(\theta) =
\frac{A_1 2^{\kappa-1}(3-\kappa)
B\left(\frac{1+\kappa}{2},2-\kappa\right)}{\kappa}.
\Label{3-8-5}
\end{align}
\end{prop}
\begin{proof}
When $A_1 \neq A_2$, (\ref{20}) and (\ref{20.1}) of Lemma \ref{l8} guarantee that
\begin{align}
\left. \frac{\,d }{\,d s}{I}^{s}_{g,\theta}
\right|_{s= \frac{1}{2}}
= (A_1-A_2)
\frac{(\kappa-1)(3-\kappa)}{4}\pi
\tan \frac{2-\kappa}{2}\pi 
B(\frac{1+\kappa}{2},2-\kappa)
s(1-s(\kappa-1))
B(s+\kappa(1-s),2-\kappa) \neq 0.
\end{align} From the concavity and the continuity of the maximized function, we have $\overline{\alpha}_1(\theta) \,> 2^{\kappa}I^{\frac{1}{2}}_{g,\theta}$. When $A_1= A_2$, we have $I^{s}_{g,\theta}=I^{1-s}_{g,\theta}$.
The relations 
\begin{align}
\overline{\alpha}_1(\theta) =
2^{\kappa}I^{\frac{1}{2}}_{g,\theta}=
\frac{A_1 2^{\kappa-1}(3-\kappa)
B\left(\frac{1+\kappa}{2},2-\kappa\right)}{\kappa}
\Label{3-8-6}
\end{align}
follow from the concavity. Since the minimums
\begin{align*}
\min_{0 \le s \le 1}
\left( s^{\frac{1}{\kappa-1}} +
(1-s)^{\frac{1}{\kappa-1}} \right)
\end{align*}
and 
\begin{align*}
\min_{0 \le s \le 1}
\left( 
\frac{(1-s(\kappa-1))B(s+\kappa(1-s), 2-\kappa)}
{1-s} +
\frac{(1-(1-s)(\kappa-1))
B((1-s)+\kappa s, 2-\kappa)}
{s} \right)
\end{align*}
are achieved at the same point $s= \frac{1}{2}$ (See Lemma \ref{l13}), the relation 
\begin{align}
&\min_{0 \le s \le 1}
\left( s^{\frac{1}{\kappa-1}} +
(1-s)^{\frac{1}{\kappa-1}} \right)^{\kappa-1}
\left( \frac{(1-s(\kappa-1))B(s+\kappa(1-s), 2-\kappa)}
{1-s} +
\frac{(1-(1-s)(\kappa-1))
B((1-s)+\kappa s, 2-\kappa)}
{s} \right)\nonumber \\
= &
2^{\kappa-1}(3-\kappa) B\left( \frac{1+\kappa}{2},
2-\kappa \right) \Label{3-8-4}
\end{align}
holds.
Thus, equation (\ref{3-8-5}) follows from (\ref{3-8-6}) and (\ref{3-8-4}).
\end{proof}

Next, we consider the case $A_2=0$.
\begin{prop}
When $1 \,< \kappa \,< 2- t_0$, 
the relations 
\begin{align}
\alpha_1(\theta)=
\overline{\alpha}_1(\theta)=
A_1 \frac{2^{\kappa}}{\kappa} \Label{3-8-7}
\end{align}
hold, where the real number 
$t_0 \in (0, \frac{1}{2})$ is uniquely
defined by (see Lemma \ref{l12})
\begin{align}
2 t_0 + t_0 (1-t_0)
\left( \psi(1+t_0) - \psi(1)\right)= 1,\Label{3-8-12}
\end{align}
where $\psi(x)$ is the D-psi function
defined by $\psi(x):= \frac{\,d }{\,d x} \log \Gamma (x)$.
\end{prop}
The number $t_0$ is enumerated by
$t_0\cong 0.432646$, as is checked by the following graph. 
\begin{figure}[hbtp]
\centering
\includegraphics{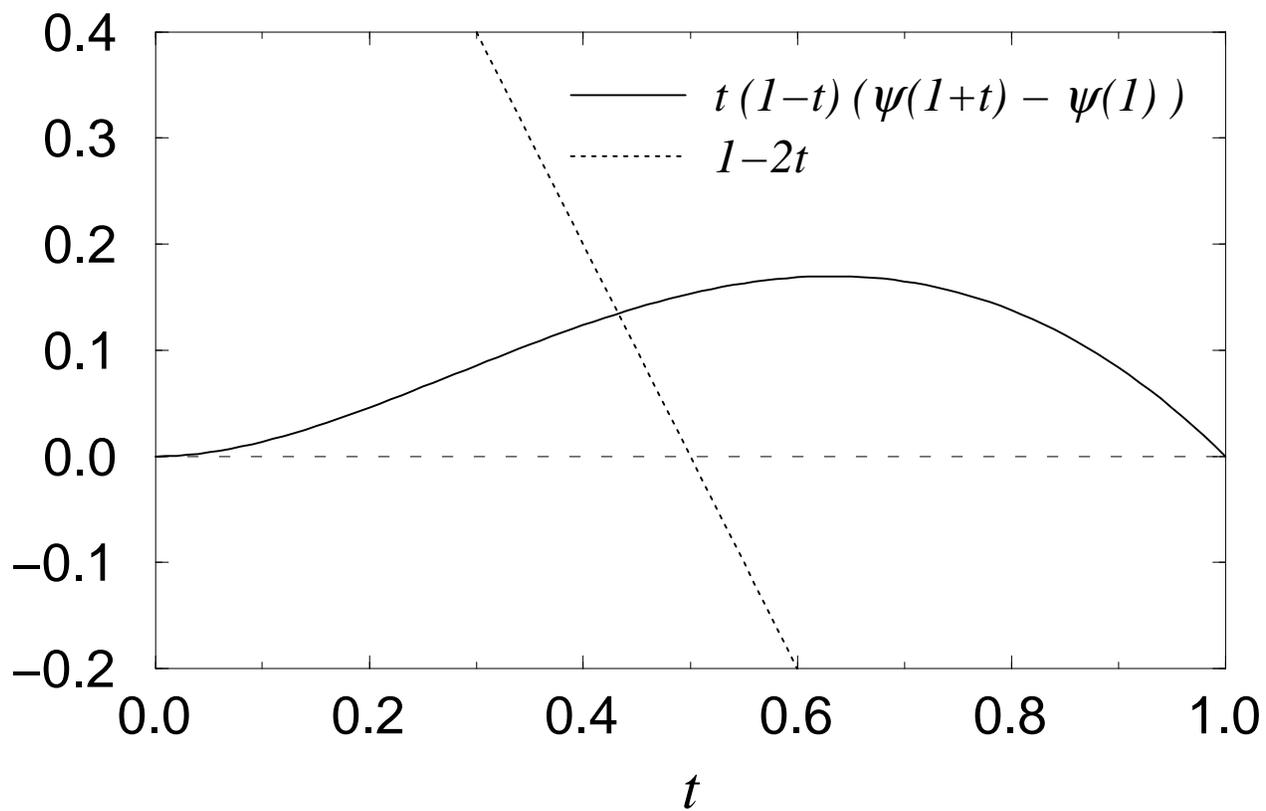}
\caption{Functions 
$t(1-t)\left( \psi(1+t) - \psi(1)\right)$ and $1-2t$}
\end{figure}

\begin{proof}
In this case, since the function $s \mapsto s(1-s(\kappa-1))B(s+\kappa(1-s), 2-\kappa)$ is concave, Lemma \ref{l12} guarantees that
\begin{align*}
\frac{\,d }{\,d s}\left. 
s(1-s(\kappa-1))B(s+\kappa(1-s), 2-\kappa)
\right)
&\ge 
\left. 
\frac{\,d }{\,d s}
s(1-s(\kappa-1))B(s+\kappa(1-s), 2-\kappa)
\right|_{s=1} \\
&=
(3-2\kappa)+(2-\kappa)(\kappa-1)
(\psi(3-\kappa)- \psi(1)) \ge 0
\end{align*}
for $s \in (0,1)$.
Thus, 
\begin{align*}
\overline{\alpha}_1 (\theta)
&= \frac{A_1 2^\kappa}{\kappa}
\max_{0 \le s \le 1}
s(1-s(\kappa-1))B(s+\kappa(1-s), 2-\kappa) \\
&= \frac{A_1 2^\kappa}{\kappa}
1 ( 1- 1 (\kappa-1))B( 1+ \kappa(1-1) ,2 - \kappa) 
= \frac{A_1 2^\kappa}{\kappa}.
\end{align*}
As in subsection \ref{e2}, we can prove that
\begin{align*}
\lim_{\epsilon \to +0}
\frac{1}{\epsilon^{\kappa}}
\beta (\vec{\underline{\theta}_\epsilon},
\theta, \epsilon) =
\frac{A_1 2^\kappa}{\kappa},
\end{align*}
which implies (\ref{3-8-7}).
\end{proof}
\subsection{The case that 
$0 \,< \kappa_1=\kappa_2 \,< 1$}\Label{e5}
As was proven by Hayashi \cite{Haya2}, when $0 \,< \kappa_1=\kappa_2 \,< 1$, the equation 
\begin{align*}
\lim_{\epsilon \to +0}
\frac{I^s( f_{\theta} \| f_{\theta+\epsilon}) }
{\epsilon^{\kappa_1}}
=
\frac{1-\kappa_1}{\kappa_1}
(A_1 s B(s + \kappa_1(1-s),1-\kappa_1)+
A_2(1-s) B(1-s + \kappa_1 s,1-\kappa_1))
\end{align*}
holds, where this convergence is uniform for $s \in (0,1)$. Letting $g(x)=x^\kappa_1$, we have $\kappa=\kappa_1$. 
\begin{align*}
&\overline{\alpha}_1(\theta)= 
2^\kappa \max_{0\le s \le 1}
I^s_{g,\theta} \\
&=
2^\kappa
\max_{0\le s \le 1}
(1-\kappa)
\frac{
A_1 s B(s + \kappa(1-s),1-\kappa)+
A_2(1-s) B(1-s + \kappa s,1-\kappa)}
{\kappa} \\
&\overline{\alpha}_2(\theta)= 
\sup_{0\,< s \,< 1}
\frac{I^s_{g,\theta}}{s(1-s)}
\left( s^{\frac{1}{\kappa-1}} 
+ (1-s)^{\frac{1}{\kappa-1}} \right)^{\kappa-1} \\
&=
\sup_{0\,< s \,< 1}
(1-\kappa)
\frac{
A_1 s B(s + \kappa(1-s),1-\kappa)+
A_2(1-s) B(1-s + \kappa s,1-\kappa)}
{\kappa s(1-s)}
\left( s^{\frac{1}{\kappa-1}} 
+ (1-s)^{\frac{1}{\kappa-1}} \right)^{\kappa-1} .
\end{align*}
\begin{prop}
If, and only if, $A_1 = A_2$, the equality
\begin{align}
\overline{\alpha}_1(\theta)
= \overline{\alpha}_2(\theta) \Label{3-8-22} 
\end{align}
holds.
In this case, the equation
\begin{align}
\overline{\alpha}_1(\theta)= 
\overline{\alpha}_2(\theta)= 
\frac{1}{\kappa} A_1 2^{\kappa}
\left( 1 - \kappa \right)
B\left(\frac{1 + \kappa}{2},1 -\kappa\right)
\Label{3-8-21} 
\end{align}
holds. 
\end{prop}
\begin{proof}
Using (\ref{21}) and (\ref{21.1}) of Lemma \ref{l8}, we obtain
\begin{align*}
\left. \frac{\,d }{\,d s}{I}^{s}_{g,\theta}
\right|_{s= \frac{1}{2}}
= (A_1-A_2)
\frac{1-\kappa}{2}\pi 
\cot \frac{1-\kappa}{2}\pi
B(\frac{1+\kappa}{2},1-\kappa).
\end{align*}
Since $
\frac{1-\kappa}{2}\pi 
\cot \frac{1-\kappa}{2}\pi
B(\frac{1+\kappa}{2},1-\kappa)\,> 0$,
Theorem \ref{thm3} yields this sufficient and necessary condition for (\ref{3-8-22}). Equation (\ref{16.3}) implies (\ref{3-8-21}).
\end{proof}

However, since the function $x \mapsto (\kappa-1)\log x$ is convex on $(0, \infty)$, the function $x \mapsto \log f(x)$ is not concave on $(a,b)$.
There does not exist an example in which Lemma \ref{thm7.2} can be applied. Thus, it is an open problem whether there exists an example such that
\begin{align*}
\overline{\alpha}_1(\theta)= 
\alpha_1(\theta)
\end{align*}
in this case, except for the case $A_1 A_2 =0$.

\begin{prop}
When $A_2=0$, $\alpha_1(\theta)$ and $ \alpha_2(\theta)$ are calculated as
\begin{align}
\alpha_1(\theta) &= 
\overline{\alpha}_1(\theta)= 
\frac{A_1 2^{\kappa}}{\kappa} \Label{3-8-10}\\
\alpha_2(\theta) &= 
\overline{\alpha}_2(\theta)= 
\frac{A_1}{\kappa} .\Label{3-8-11}
\end{align}
\end{prop}
\begin{proof}
Since the function $s \mapsto s B(s+\kappa(1-s), 1-\kappa)$ is concave, we have
\begin{align}
\frac{\,d }{\,d s}
s B(s+\kappa(1-s), 1-\kappa)
&\ge
\left.\frac{\,d }{\,d s}
s B(s+\kappa(1-s), 1-\kappa)\right|_{s=1} \nonumber\\
&=
1+ (1-\kappa) ( \psi(1)-\psi(2-\kappa))\nonumber \\
& \ge 
1+ (1-\kappa) ( \psi(1)-\psi(2)) \Label{3-8-8}\\
&= 1- (1-\kappa) = \kappa \,> 0,\Label{3-8-9}
\end{align}
where inequality (\ref{3-8-8}) holds because $\psi(x)$ is monotonically increasing in $x \in (0, \infty)$, and the first equation of (\ref{3-8-9}) follows from the formula $\psi(x+1)= \psi(x) + \frac{1}{x}$. Thus,
\begin{align*}
\overline{\alpha}_1(\theta)
&=
\frac{A_1 2^{\kappa}(1-\kappa)}{\kappa}
\max_{0 \le s \le 1}
s B(s+\kappa(1-s), 1-\kappa) \\
&= \frac{A_1 2^{\kappa}(1-\kappa)}{\kappa}
B(1, 1-\kappa) 
= \frac{A_1 2^{\kappa}}{\kappa}.
\end{align*}
Since, as in subsection \ref{e4}, we can check that the estimators $\{ \underline{\theta}_\epsilon\} _{\epsilon \,> 0}$ achieve the bound $\frac{A_1 2^{\kappa}}{\kappa}$, equations (\ref{3-8-10}) hold.

The other upper bound $\overline{\alpha}_2(\theta)$ is calculated as 
\begin{align*}
\overline{\alpha}_2(\theta)
=
\frac{A_1 (1-\kappa)}{\kappa}
\max_{0 \le s \le 1}
B(s+\kappa(1-s), 1-\kappa) 
\left( \left( \frac{1-s}{s}\right)^{\frac{1}{1-\kappa}}
+1\right)^{-(1-\kappa)}.
\end{align*}
Note that the beta function $B(x,y)$ is monotonically decreasing for $x,y \,> 0$. Since both 
$\max_{0\le s \le 1}
B(s+\kappa(1-s), 1-\kappa) $
and 
$\max_{0\le s \le 1}
\left( 
\left( \frac{1-s}{s}\right)^{\frac{1}{1-\kappa}}
+1\right)^{-(1-\kappa)}$
are achieved at the same point, $s=1$, we have 
\begin{align*}
\overline{\alpha}_2(\theta)
= \frac{A_1(1-\kappa)}{\kappa}B(1,1-\kappa)
= \frac{A_1}{\kappa}.
\end{align*}
This bound is achieved by the estimator $\vec{\underline{\theta}}$ because 
\begin{align*}
\lim_{\epsilon \to +0}
\frac{1}{\epsilon^{\kappa}}
\beta(\vec{\underline{\theta}},\theta,\epsilon)
=
\lim_{\epsilon \to +0}
\frac{-1}{\epsilon^{\kappa}}
\log \int_{a+\epsilon}^b
f(x) \,d x=
\frac{A_1}{\kappa}.
\end{align*}
Therefore, we have (\ref{3-8-11}).
\end{proof}

\section{Conclusion}
We have discussed large deviation theories under a more general setting. The two quantities $\overline{\alpha}_1(\theta)$ and $\overline{\alpha}_2(\theta)$ do 
not necessarily coincide. In a non-regular case, it is clear that the order of limits is crucial. In the future, 
such phenomena deserve study from another viewpoint.

Nagaoka \cite{Naga2} initiated a discussion of two kinds 
of large deviation bounds, as in this paper, in a quantum setting, and Hayashi \cite{Haya} discussed these in more 
depth. The two kinds of large deviation bounds do not necessarily coincide in a quantum setting. However, the reason for this difference in a quantum setting differs from that 
for a non-regular setting. Gaining an understanding of 
these differences from a unified viewpoint remains a goal for the future.

\appendix
\section{Lemmas concerning the beta and D-psi functions}
In this section, we prove some formulas concerning the beta and D-psi functions used in Section \ref{s6}.
\begin{lem} \Label{l8}
When $1 \,< \kappa \,< 2$, we have 
\begin{align}
\left. \frac{\,d }{\,d s }
s(1-s(\kappa-1))B(s+\kappa(1-s),2-\kappa)
\right|_{s=1/2} &=
\frac{(\kappa-1)(3-\kappa)}{4}\pi
\tan \frac{2-\kappa}{2}\pi B(\frac{1+\kappa}{2},2-\kappa)
\Label{20} \,> 0\\
\left. \frac{\,d }{\,d s }
(1-s)(1-(1-s)(\kappa-1))B((1-s)+\kappa s,2-\kappa)
\right|_{s=1/2} &
=
- 
\frac{(\kappa-1)(3-\kappa)}{4}\pi
\tan \frac{2-\kappa}{2}\pi B(\frac{1+\kappa}{2},2-\kappa)
\,< 0
\Label{20.1} 
\end{align}
When $0 \,< \kappa \,< 1$, the equations
\begin{align}
\left. \frac{\,d }{\,d s }
s B(s+\kappa(1-s),1-\kappa)
\right|_{s=1/2} &=
\frac{1-\kappa}{2}\pi 
\cot \frac{1-\kappa}{2}\pi
B(\frac{1+\kappa}{2},1-\kappa)\,> 0
\Label{21} \\
\left. \frac{\,d }{\,d s }
(1-s) B((1-s)+\kappa s,1-\kappa)
\right|_{s=1/2} &=
- \frac{1-\kappa}{2}\pi 
\cot \frac{1-\kappa}{2}\pi
B(\frac{1+\kappa}{2},1-\kappa) \,< 0
\Label{21.1}
\end{align}
hold.
\end{lem}
\begin{proof}
Using the function 
$\psi(x):= \frac{\,d }{\,d x} \log \Gamma (x)$,
we can calculate
\begin{align}
& \left. \frac{\,d }{\,d s }
s(1-s(\kappa-1))B(s+\kappa(1-s),2-\kappa)
\right|_{s=1/2} \nonumber\\
&=
\left((2-\kappa)
+ \frac{(3-\kappa)(1-\kappa)}{4}
\left( \psi(\frac{1+\kappa}{2})
- \psi(\frac{5-\kappa}{2})\right)
\right)
B(\frac{1+\kappa}{2},2-\kappa) \nonumber\\
&=
\left((2-\kappa)
+ \frac{(3-\kappa)(1-\kappa)}{4}
\left( \psi(\frac{\kappa-1}{2}) + \frac{2}{\kappa-1}
- \psi(\frac{3-\kappa}{2}) - \frac{2}{3-\kappa}
\right)
\right)
B(\frac{1+\kappa}{2},2-\kappa) \Label{23}\\
&=
\left((2-\kappa)
+ \frac{(3-\kappa)(1-\kappa)}{4}
\left( 
\pi \cot \pi \frac{3-\kappa}{2}
+ \frac{8-4\kappa}{(\kappa-1)(3-\kappa)}
\right)
\right)
B(\frac{1+\kappa}{2},2-\kappa) \Label{24}\\
&=
\frac{(\kappa-1)(3-\kappa)}{4}\pi
\tan \frac{2-\kappa}{2}\pi B(\frac{1+\kappa}{2},2-\kappa)
\Label{25}
\end{align}
where we use the formula 
$\psi(x+1)= \frac{1}{x}+\psi(x)$ in (\ref{23}),
the formula 
$\psi(1-x) - \psi(x) = \pi \cot \pi x$ in (\ref{24}),
and the formula 
$\cot (\frac{\pi}{2} +x)= - \tan x $ in (\ref{25}).
We obtain (\ref{20}).
Similarly, we can prove (\ref{20.1}).

Next, we prove (\ref{21}). We can calculate
\begin{align}
&\left. \frac{\,d }{\,d s }
s B(s+\kappa(1-s),1-\kappa)
\right|_{s=1/2} \nonumber\\
& =
\left( 1 + \frac{1-\kappa}{2}
\left( \psi(\frac{1+\kappa}{2}) 
- \psi(\frac{3-\kappa}{2}) \right)
\right)
B(\frac{1+\kappa}{2}, 1-\kappa) \nonumber \\
& = \left( 1 + \frac{1-\kappa}{2}
\left( \psi(\frac{1+\kappa}{2}) 
- \psi(\frac{3-\kappa}{2}) -
\frac{2}{1-\kappa}
\right)
\right)
B(\frac{1+\kappa}{2}, 1-\kappa) \Label{27}\\
& = 
\frac{1-\kappa}{2}
\pi \cot \frac{1-\kappa}{2}\pi 
B(\frac{1+\kappa}{2}, 1-\kappa) \Label{26},
\end{align}
where (\ref{27}) follows from the formula $\psi(x+1)= \frac{1}{x}+\psi(x)$ and (\ref{26}) follows from the formula $\psi(1-x) - \psi(x) = \pi \cot \pi x$.
Similarly, we obtain (\ref{21.1}).
\end{proof}

\begin{lem}\Label{l12}
Assume that $1 \,< \kappa \,<2$. There uniquely exists $t_0 \in (0, \frac{1}{2})$ satisfying (\ref{3-8-12}), and the inequality
\begin{align}
(3-2\kappa)+(2-\kappa)(\kappa-1)
(\psi(3-\kappa)- \psi(1)) \ge 0 \Label{3-8-13}
\end{align}
holds if, and only if, $1\,< \kappa \le 2- t_0$.
\end{lem}
\begin{proof}
In the following, this lemma is proven by replacing $\kappa$ with $2-t$.
Define the function $h(t):=
2t-1 +t(1-t) (\psi( 1+t) -\psi(1))$.
When $t \,< \frac{1}{2}$,
we have 
\begin{align*}
h'(t)=
2+ (1-2t)(\psi(1+t) -\psi(1))
+(t-2t^2)\psi'(1+t) \,> 0
\end{align*}
because 
$\psi'(x) \ge 0$ for $x \,> 0$.
Therefore, 
$h(t)$ is strictly monotonically increasing in 
$(0,\frac{1}{2})$.
Since
$h(0)=-1\,<0, h(\frac{1}{2})=\frac{1}{4}
\left(\psi(\frac{3}{2})-\psi(1)\right) \,>0$,
there uniquely exists the number $t_0 \in (0 ,\frac{1}{2})$ satisfying (\ref{3-8-12}).
Also, in this case, the inequality $h(t) \ge 0$ holds if $t \ge t_0$.

Next, we consider case $t \ge \frac{1}{2}$.
Given the relations,
\begin{align*}
2t-1 \ge 0 ,\quad
t(1-t) \ge 0, \quad
\psi(1+t) -\psi(1) \ge 0, 
\end{align*}
$h(t) \ge 0$.
\end{proof}

\begin{lem}\Label{l13}
The minimum
\begin{align*}
\min_{0 \le s \le 1}
\left(
\frac{(1-s(\kappa-1))B(s+\kappa(1-s), 2-\kappa)}
{1-s} +
\frac{(1-(1-s)(\kappa-1))
B((1-s)+\kappa s, 2-\kappa)}
{s} \right)
\end{align*}
is attained at $s = \frac{1}{2}$.
\end{lem}
\begin{proof}
Since the minimized function is invariant for the replacement $s \mapsto 1-s$, it is sufficient to show its concavity. Since 
$\frac{\,d^2}{\,d x^2}
e^{h(x)} =( h''(x)+(h('x))^2)e^{h(x)}$,
we can show its concavity by proving the concavity of the function
$s \mapsto
\log \frac{1-s(\kappa-1)}{1-s}
B(s+\kappa(1-s), 2-\kappa)$.
We can evaluate 
\begin{align*}
\frac{\,d^2}{\,d s^2}
\log B(s+\kappa(1-s), 2-\kappa)
= (1-\kappa)^2
\left( \psi'( s+\kappa(1-s))
-\psi'( s+\kappa(1-s)+2 -\kappa)\right)
\,>0
\end{align*}
because $\psi'(x)$ is monotonically decreasing for $x \,>0$.
Also, we have 
\begin{align*}
\frac{\,d^2}{\,d s^2}
\log\frac{1-s(\kappa-1)}{1-s}
= \frac{(2-\kappa) (\kappa-2 (\kappa-1)s)}
{(1+ \kappa s + (\kappa-1)s^2)^2} \,> 0
\end{align*}
because $\kappa \,> 2(\kappa -1)$. 
The proof is now complete.
\end{proof}
\section{Concave function}
\begin{lem}\Label{concave}
When a concave function $f \ge 0$ is defined in $(0,1)$,
\begin{align*}
\inf_{x \ge0}
sx +(1-s)\sup_{0 \,< t \,< 1}
\frac{-tx +f(t)}{1-t}
=
\inf_{x \,>0}\sup_{0 \,< t \,< 1}
\frac{(s-t)x +(1-s)f(t)}{1-t}
= f(s).
\end{align*}
\end{lem}
\begin{proof}
Substituting $s$ into $t$ we have
\begin{align*}
f(s) \le
\sup_{0 \,< t \,< 1}
\frac{(s-t)x +(1-s)f(t)}{1-t}.
\end{align*}
Taking the infimum $\inf_{x \,>0}$, we obtain 
\begin{align*}
f(s) \le \inf_{x \,>0}\sup_{0 \,< t \,< 1}
\frac{(s-t)x +(1-s)f(t)}{1-t}.
\end{align*}

Next, we proceed to the opposite inequality. From the concavity of $f$, we can define the upper derivative $\overline{f'}$ and the lower derivative $\underline{f'}$ as
\begin{align*}
\overline{f'}(s):= \lim_{\epsilon to +0}
\frac{f(s)- f(s-\epsilon)}
{\epsilon} ,\quad 
\underline{f'}(s):= \lim_{\epsilon to +0}
\frac{f(s+\epsilon)- f(s)}
{\epsilon} 
\end{align*}
Since the concavity guarantees that
\begin{align*}
\frac{f(s+\epsilon)- f(s)}
{\epsilon} 
\le \underline{f'}(s) \le \overline{f'}(s)\le 
\frac{f(s)- f(s-\epsilon)}
{\epsilon} , \quad \forall \epsilon\,> 0 ,
\end{align*}
we obtain
\begin{align*}
&\frac{(s-s)(f(s)+(1-s)\overline{f'}(s))
+(1-s)f(s)}{1-s}
- 
\frac{(s-t)(f(s)+(1-s)\overline{f'}(s))
+(1-s)f(t)}{1-t} \\
=& \frac{1-s}{1-t}
\left( f(s)-f(t) +(t-s)\overline{f'}(s)\right)
\ge 0,\quad \forall s, \forall t \in (0,1).
\end{align*}
Therefore,
\begin{align*}
f(s) \ge \sup_{0\,< t \,< 1}
\frac{(s-t)(f(s)+(1-s)\overline{f'}(s))
+(1-s)f(t)}{1-t}
\ge 
\inf_{x \,> 0}
\sup_{0\,< t \,< 1}
\frac{(s-t)(f(s)+(1-s)\overline{f'}(s))
+(1-s)f(t)}{1-t}.
\end{align*}
The proof is now complete.
\end{proof}

\section{Other lemmas}\Label{A}
\begin{lem}\Label{ap1}
When $g$ is strictly monotonically decreasing and continuous and satisfies that $g(0)=0$, there exists $\kappa \,> 0$ such that
\begin{align}
x^\kappa = \lim_{\epsilon \to +0}
\frac{g(x\epsilon)}{g(\epsilon)}, \quad
x \,> 0. \Label{82}
\end{align}
\end{lem}
\begin{proof}
Let $h(x)$ be the RHS of (\ref{82}).
Since 
\begin{align*}
\lim_{\epsilon \to +0}
\frac{g(xy \epsilon)}{g(\epsilon)}
=
\lim_{\epsilon \to +0}
\frac{g(xy \epsilon)}{g(y\epsilon)}
\lim_{\epsilon \to +0}
\frac{g(y \epsilon)}{g(\epsilon)},
\end{align*}
$h(xy)=h(x)h(y)$.
Thus, there exists $\kappa \,> 0$ satisfying (\ref{82}).
\end{proof}
\begin{lem}\Label{l11}
For any $\delta \,>0$, we have
\begin{align}
\lim_{\epsilon \to 0}
\frac{1}{\epsilon^2 \log \epsilon}
\int_0^\delta \exp \left( -\epsilon
\frac{1}{x}\right)(x+\epsilon) - x \,d x
&= \frac{1}{2} \Label{28-8}\\
\lim_{\epsilon \to 0}
\frac{1}{ \epsilon^2 \log \epsilon}\
\int_\epsilon^\delta \exp \left( +\epsilon
\frac{1}{x}\right)
(x-\epsilon) - x \,d x
&= \frac{1}{2} \Label{28-9}.
\end{align}
\end{lem}
\begin{proof}
we can calculate
\begin{align*}
\int_\epsilon^\delta \exp \left( +\epsilon
\frac{1}{x}\right)
(x-\epsilon) - x \,d x
=\int_\epsilon^\delta 
\sum_{n=1}^\infty - 
\frac{n \epsilon^{n+1}}{(n+1)! x^n} \,d x 
=
-\frac{\epsilon^2}{2}(\log \delta - \log \epsilon)
+\sum_{n=2}^\infty
\frac{n(n-1)}{(n+1)!}
\left( \frac{\epsilon^{n+1}}{\delta^{n-1}}
-\epsilon^2\right).
\end{align*}
Since $\frac{\epsilon^2}{\epsilon^2 \log \epsilon }
\to 0$, we obtain (\ref{28-9}). Similarly, we can calculate
\begin{align*}
\int_\epsilon^\delta \exp \left( -\epsilon
\frac{1}{x}\right)(x+\epsilon) - x \,d x 
=
-\frac{\epsilon^2}{2}(\log \delta - \log \epsilon)
+\sum_{n=2}^\infty
\frac{n(n-1)(-1)^n}{(n+1)!}
\left( \frac{\epsilon^{n+1}}{\delta^{n-1}}
-\epsilon^2\right).
\end{align*}
Since
\begin{align*}
0 \le \int_0^\epsilon \exp \left( -\epsilon
\frac{1}{x}\right)(x+\epsilon) \,d x 
& \le\int_0^\epsilon(x+\epsilon) \,d x 
= \frac{3}{2}\epsilon^2 \\
\int_0^\epsilon x \,d x &= \frac{1}{2}\epsilon^2,
\end{align*}
we obtain (\ref{28-8}).
\end{proof}

\begin{lem}\Label{L38}
For $0 \,< s \,< \frac{1}{2}$, the inequalities
\begin{align*}
2s I^{\frac{1}{2}}(p\|q)
\le I^{s}(p\|q)
\le 2(1-s) I^{\frac{1}{2}}(p\|q)
\end{align*}
hold.
\end{lem}
\begin{proof}
Since $\left(\frac{1}{2s}\right)^{-1}+
\left(\frac{1}{1-2s}\right)^{-1}=1$,
the H\"{o}lder inequality guarantees that
\begin{align*}
\int_{\Omega}
p^{s}(\omega) q^{1-s}(\omega) \,d \omega
&=
\int_{\Omega}
\left(p^{s}(\omega) q^{s}(\omega) \right)
\left(q^{1-2s}(\omega)\right) \,d \omega \\
& \le
\left(\int_{\Omega}
\left(p^{s}(\omega) q^{s}(\omega) \right)^{\frac{1}{2s}}
\,d \omega \right)^{2s}\cdot
\left(
\int_{\Omega}
\left(q^{1-2s}(\omega)\right)^{\frac{1}{1-2s}}
 \,d \omega 
\right)^{1-2s}\\
& = \left(\int_{\Omega}
p^{\frac{1}{2}}(\omega) q^{\frac{1}{2}}(\omega) 
\,d \omega \right)^{2s}.
\end{align*}
Thus, we obtain 
\begin{align*}
I^{s}(p\|q) \ge 2s I^{\frac{1}{2}}(p\|q).
\end{align*}
Similarly, since 
$\left( 2(1-s)\right)^{-1}+
\left( \frac{2(1-s)}{1-2s}\right)^{-1}=1$,
we can apply the H\"{o}lder inequality as
\begin{align*}
\int_{\Omega}
p^{\frac{1}{2}}(\omega) q^{\frac{1}{2}}(\omega) 
\,d \omega 
&=\int_{\Omega}
\left(p^{s}(\omega) q^{1-s}(\omega) 
\right)^{\frac{1}{2(1-s)}}
p^{\frac{1-2s}{2(1-s)}}(\omega) 
\,d \omega \\
&\le
\left(\int_{\Omega}
\left(p^{s}(\omega) q^{1-s}(\omega) 
\right)^{\frac{1}{2(1-s)}\cdot 2(1-s)}
\,d \omega \right)^{\frac{1}{2(1-s)}}\cdot
\left(\int_{\Omega}
p^{\frac{1-2s}{2(1-s)}\cdot\frac{2(1-s)}{1-2s}}(\omega) 
\,d \omega \right)^{\frac{1-2s}{2(1-s)}} \\
&=
\left(\int_{\Omega}
p^{s}(\omega) q^{1-s}(\omega) \,d \omega
\right)^{\frac{1}{2(1-s)}},
\end{align*}
which implies that
\begin{align*}
I^{\frac{1}{2}}(p\|q) \ge \frac{1}{2(1-s)}
I^{s}(p\|q).
\end{align*}
\end{proof}

\end{document}